\documentclass[10pt]{article}

\usepackage{latexsym}
\usepackage{amsfonts}
\usepackage{amsmath}
\usepackage{amssymb}
\newtheorem{thm}{Theorem}
\newtheorem{lem}[thm]{Lemma}

\newtheorem{prop}[thm]{Proposition}

\def\eps{\varepsilon}

\def\qed{\hfill $\vcenter{\hrule height .3mm
\hbox {\vrule width .3mm height 2.1mm \kern 2mm \vrule width .3mm
height 2.1mm} \hrule height .3mm}$ \bigskip}
\def\P{\mathbb{P}}

\def\EE{\mathbb{E}}
\def\cE{\mathcal{E}}

\def\R{\mathbb{R}}
\def\RR{\mathbb{R}}

\title{\textsf{Pointwise Estimates for Marginals of Convex Bodies}}
\author{R. Eldan \ and \ B. Klartag\thanks{Supported by the Clay Mathematics Institute and by  NSF grant
\#DMS-0456590.}}
\date{}
\begin{document}
\maketitle

\begin{abstract}
We prove a pointwise version of the multi-dimensional central limit theorem for convex bodies.
Namely, let $\mu$ be an isotropic, log-concave
 probability
measure on $\R^n$. For a typical subspace $E \subset \R^n$
of dimension $n^c$, consider the probability density of the projection of
$\mu$ onto $E$.
We show that the ratio between this probability density and the standard gaussian
density in $E$ is very close to $1$ in large parts of $E$.
Here $c > 0$ is a universal constant.
This complements a recent result by the second named author,
where the total-variation metric between the densities was considered.
\end{abstract}

\section{Introduction}

Suppose $X$ is a random vector in $\R^n$ that is distributed uniformly in some convex
set $K \subset \R^n$. For a subspace $E \subset \R^n$ we denote by $Proj_E$ the
orthogonal projection operator onto $E$ in $\R^n$.
The central limit theorem for convex bodies \cite{K1, K2} asserts
that there exists a subspace $E \subset \R^n$, with $\dim(E) > n^c$, such that
the random vector $Proj_E(X)$ is approximately gaussian, in the total variation sense.
This means that for a certain gaussian random vector $\Gamma$ in the subspace $E$,
\begin{equation} \sup_{A \subseteq E} \, \left|   \,  \P \left \{  Proj_E(X) \in A \right \} \, - \, \P \left \{ \Gamma \in A \right \} \, \right| \leq \frac{C}{n^c}, \label{eq_1148} \end{equation}
where the supremum runs over all measurable subsets $A \subseteq E$. Here, and throughout this note,
the letters $c, C, c_1, C_2, c^{\prime}, \tilde{C}$ etc. denote some positive universal constants, whose
value may change from one appearance to the next.

\medskip The total variation estimate (\ref{eq_1148}) implies
that the density of $Proj_E(X)$ is close to the density of $\Gamma$ in the $L^1$-norm.
In this note we observe that a stronger conclusion is within reach: One may deduce
that the ratio between the density of $Proj_E(X)$ and the density of $\Gamma$ deviates
from $1$ by no more than $C n^{-c}$, in the significant parts of the subspace $E$.

\medskip We need some notation. Write $| \cdot |$ for the standard Euclidean norm in $\RR^n$.
A random vector $Z$ in $\RR^n$ is isotropic if
the following normalization holds:
\begin{equation}
 \EE Z = 0, \ \ \ \ \ Cov(Z) = Id
\label{eq_1138}
 \end{equation}
where $Cov(Z)$ stands for the covariance matrix of $Z$, and $Id$ is the identity matrix.
The grassman manifold
$G_{n,\ell}$ of all $\ell$-dimensional subspaces of $\RR^n$ carries a unique rotationally-invariant
probability measure $\mu_{n,\ell}$.
Whenever we say that $E$ is a random $\ell$-dimensional subspace
in $\RR^n$,  we relate to the above probability measure $\mu_{n,\ell}$.
Under the additional assumption that the random vector $X$ is isotropic, the subspace $E$ for which $Proj_E(X)$
is approximately gaussian may be chosen at random \cite{K1, K2}.

\medskip
A function $f: \RR^n \rightarrow [0, \infty)$ is log-concave if $\log f: \RR^n \rightarrow
[-\infty, \infty)$ is a concave function. The characteristic function of a convex
set is log-concave. Throughout the entire discussion, the requirement that $X$ is
distributed uniformly in a convex body could have been relaxed to the weaker
condition, that $X$ has a log-concave density. Our main result in this paper reads as follows:
\begin{thm}
Let $X$ be an isotropic random vector in $\R^n$ with a log-concave density.
Let $1 \leq \ell \leq n^{c_1}$ be an integer. Then there exists a subset $\cE \subseteq G_{n,\ell}$ with
$\mu_{n,\ell}(\cE) \geq 1 - C \exp( -n^{c_2} )$ such that for any $E \in \cE$, the following
holds: Denote  by $f_E$ the density of the random vector $Proj_E(X)$. Then,
\begin{equation}
\left | \frac {f_E(x)}{\gamma(x)} - 1 \right | \leq \frac{C}{n^{c_3}}
\label{eq_1241}
\end{equation}
for all $x \in E$ with $|x| \leq n^{c_4}$. Here, $\gamma(x) = (2 \pi)^{-\ell/2} \exp( -|x|^2  / 2)$
is the standard gaussian density in $E$, and $C, c_1,c_2,c_3,c_4 > 0$ are universal constants.
\label{thm_main}
\end{thm}

Note that almost the entire mass of a standard $\ell$-dimensional gaussian distribution
is contained in a ball of radius $10 \sqrt{\ell}$ about the origin. Therefore, (\ref{eq_1241})
easily implies the total variation bound mentioned above.
The history of the central limit theorem for convex bodies goes
back to the conjectures and results of Brehm and Voigt \cite{BV} and
Anttila, Ball and Perissianki \cite{ABP}, see \cite{K1} and references therein.
The case $\ell = 1$ of Theorem \ref{thm_main} was proved in \cite{K2} using
the moderate deviation estimates of Sodin \cite{Sod}.
The generalization to higher dimensions is the main contribution  of the present
paper. See also \cite{bb} and \cite{ABBP}.

\medskip The basic idea of the proof of Theorem \ref{thm_main}
is the following: It is shown in \cite{K2}, using concentration techniques,
that the density of $Proj_E(X + Y)$ is pointwise approximately radial, where $Y$ is an independent small gaussian random vector. It is furthermore proved that this density is concentrated in a thin spherical shell.
We combine these facts to
deduce, in Section \ref{sec2}, that the density of $Proj_E(X + Y)$ is not only radial, but
in fact very close to the gaussian density in $E$.
Then, in Section \ref{sec3}, we show that the addition of the gaussian random vector $Y$
is not required. That is, we prove that when a log-concave density convolved with a small gaussian is almost gaussian -- then the original
density is also approximately gaussian. This completes the sketch of the proof.

\medskip
\emph{Acknowledgements}. The first named author would like to
express his sincere gratitude to his supervisor, Prof. Vitali Milman
who introduced him to the subject, guided him and encouraged him to
write this note. We would also like to thank Sasha Sodin and Prof.
Vitali Milman for reviewing a preliminary version of this note.

\section{Convolved marginals are Gaussian}
\label{sec2}

For a dimension $n$ and $v > 0$ we write
\begin{equation} \gamma_{n}[v](x) = \frac{1}{(2 \pi v)^{n/2}} \exp
\left( -\frac{|x|^2}{2v} \right) \ \ \ \ \ \ \ \ \ \ \ (x \in \RR^n).
\end{equation}
 That is, $\gamma_{n}[v]$ is the density of a gaussian random
vector in $\RR^n$ with mean zero and covariance  matrix $v
Id$.  Let $X$ be an isotropic random vector with a log-concave density in $\R^n$, and let $Y$ be an
independent
gaussian random vector in $\RR^n$ whose density is $\gamma_n[ n^{-\alpha}]$,
for a parameter $\alpha$ to be specified later on.
Denote by $f_{X+Y}$ the density of the random vector $X + Y$.
Our first step is to show that the density of the projection of $X + Y$ onto a typical
subspace is pointwise approximately gaussian.

\medskip
We follows the notation of \cite{K2}. For an integrable function $f: \RR^n
\rightarrow [0, \infty)$, a subspace $E \subseteq \RR^n$ and a
point $x \in E$ we write
\begin{equation} \pi_E(f)(x) = \int_{x + E^{\perp}} f(y) dy, \end{equation}
where $x + E^{\perp}$ is the affine subspace orthogonal to $E$ that passes through
the point $x$. In other words, $\pi_E(f): E \rightarrow [0, \infty)$ is the marginal of
$f$ onto $E$. The group of all orthogonal transformations of
determinant one in $\RR^n$ is denoted by $SO(n)$. Fix a dimension
$\ell$ and
a subspace $E_0 \subset \RR^n$ with $\dim(E_0) = \ell$.
For $x_0 \in E_0$ and a rotation $U \in SO(n)$, set
\begin{equation}
M_{f,E_0,x_0}(U) = \log \pi_{E_0} (f \circ U)(x_0).
\end{equation}
Define
\begin{equation}
M(|x_0|) = \int_{SO(n)} M_{f_{X+Y}, E_0, x_0} (U) d \mu_n(U),
\end{equation}
where $\mu_n$ stands for the unique rotationally-invariant Haar probability
measure on $SO(n)$. Note that $M(|x_0|)$ is independent of the
direction of $x_0$, so it is well defined. We learned in
\cite{K2} that the function $U \mapsto M_{f_{X+Y},E_0,x_0}(U)$ is highly concentrated
with respect to $U$ in the special orthogonal group $SO(n)$, around
its mean value $M(|x_0|)$. This implies that the function
$\pi_E(f_{X+Y})$ is almost spherically symmetric, for a typical
subspace $E$. This information is contained in our next Lemma, which
is equivalent to \cite[Lemma 3.3]{K2}.

\begin{lem} \label{concent}
Let $1 \leq \ell \leq n$ be integers, let $0 < \alpha < 10^5$ and denote $\lambda = \frac{1}{5 \alpha + 20}$.
Assume that $\ell \leq n^{\lambda}$. Suppose that $X$ is an isotropic random
vector with a log-concave density and that $Y$ is an independent
random vector with density $\gamma_n[n^{- \alpha \lambda}]$. Denote
the density of $X+Y$ by $f_{X+Y}$.

\medskip Let $E \in G_{n,\ell}$ be a
random subspace. Then, with probability greater than $1-C
e^{-c n^{1/10}}$ of selecting $E$, we have
\begin{equation}
\left | \log \pi_E(f_{X+Y})(x)  - M(|x|) \right | \leq
Cn^{-\lambda}, \label{eq_1238}
\end{equation}
for all $x \in E$ with $|x| \leq 5 n^{\lambda/2}$. Here $c, C > 0$ are universal constants.
\end{lem}
\textbf{Sketch of Proof:} We have to follow the proof of Lemma 3.3
in \cite{K2}, choosing for instance, $u = \frac{9}{10}$, $\lambda =
\frac{1}{5 \alpha + 20}$, $k = n^{\lambda}$ and $\eta = 1$.
Throughout the argument in \cite{K2}, it was assumed that the
dimension of the subspace is exactly $k = n^{\lambda}$, while in the
present version of the statement, note that it could possibly be
smaller, i.e., $\ell \leq k$ (note also that here, $k$ need not be an
integer). We re-run the proofs of Lemmas 2.7, 2.8, 3.1 and 3.3 from
\cite{K2}, allowing the dimension of the subspace we are working
with to be smaller than $k$,
 noting that the reduction of the dimension always acts in our benefit. \\
We refer the reader to the original argument in the proof of
Lemma 3.3 in \cite{K2} for more details. \hfill $\square$

\medskip Our main goal in this section is to show that $M(|x|)$ behaves
approximately like $\log \gamma_n[1+n^{-\alpha \lambda}](x)$. Once we prove this, it
would follow from the above lemma that the density of $X+Y$ is
pointwise approximately gaussian. Next we explain why no serious
harm is made if we take the logarithm outside the integral in the definition of $M(|x|)$.
Denote, for $x \in E_0$,
\begin{equation}
 \tilde{M}(|x|) = \int_{SO(n)} \pi_{E_0} (f_{X+Y} \circ U)(x)  d \mu_n(U).
 \label{eq_942}
 \end{equation}

\begin{lem} Under the notation and assumptions of Lemma \ref{concent}, for $|x| \leq 5 n^{\lambda / 2}$
we have
\begin{equation}
0 \leq \log \tilde{M}(|x|) - M(|x|) \leq \frac{C}{n^{1/5}},
\label{eq_926}
  \end{equation}
where $C > 0$ is a universal constant.
\label{nolog}
\end{lem}

\emph{Proof:} Recall that $E_0 \subset \RR^n$ is some fixed $\ell$-dimensional subspace.
Fix $x_0 \in E_0$ with $|x_0| \leq 5 n^{\lambda / 2}$. Lemma 3.1 of \cite{K2} states that
for any   $U_1, U_2 \in SO(n)$,
\begin{equation}
 \left| M_{f_{X+Y}, E_0, x_0}(U_1) - M_{f_{X+Y}, E_0, x_0}(U_2)  \right| \leq
C n^{\lambda (2 \alpha + 2)}  \cdot d(U_1, U_2),
\label{eq_720}
\end{equation}
where $d(U_1,U_2)$ stands for the geodesic distance between $U_1$ and $U_2$ in $SO(n)$.
As we mentioned before, Lemma 3.1 is proved in \cite{K2} under the assumption that
the dimension of the subspace $E_0$ is exactly $n^{\lambda}$.
In our case, the dimension $\ell$ might be smaller than $n^{\lambda}$, but
a direct inspection of the proofs in \cite{K2} reveals that the reduction
of the dimension can only improve the estimates. Hence (\ref{eq_720}) holds true.

\medskip
We apply the Gromov-Milman concentration inequality on $SO(n)$, quoted as Proposition
3.2 in \cite{K2}, and conclude from (\ref{eq_720}) that for any $\eps > 0$,
\begin{equation} \mu_n \left \{ U \in SO(n) ; \left| M_{f_{X+Y}, E_0, x_0}(U) - M(|x_0|) \right| \geq \eps \right \}
\leq \bar{C} \exp \left(- \bar{c} n \eps^2 / L^2 \right),  \label{eq_729}
\end{equation}
with $L = C n^{\lambda (2 \alpha + 2)}$. That is, the distribution
of $$ F(U) = \frac{\sqrt{n}}{L} \left( M_{f_{X+Y}, E_0, x_0}(U) -
M(|x_0|) \right) \ \ \ \ \ \ (U \in SO(n)) $$ on $SO(n)$ has a
subgaussian tail. Note also that $\int_{SO(n)} F(U) d \mu_n(U) = 0$. A standard  computation shows for any $p \geq 1$,
\begin{equation} \int_{SO(n)} F^p(U) d \mu_n(U)  \leq \left( C^{\prime} \sqrt{p} \right)^p, \end{equation}
where $C'$ is a universal constant. Hence, for any $0 < t < 1$,
\begin{eqnarray}
\label{eq_933}
\lefteqn{ \int_{SO(n)}  \exp \left( t   F(U) \right) d\mu_n(U) \leq
 1 + t \int_{SO(n)} F(U) d \mu_n(U)   + \sum_{i=2}^{\infty} \left( C^{\prime} \sqrt{i} \right)^i \frac{t^i}{i!} } \\
 & \leq  & 1 + \sum_{i=2}^{\infty}
 \frac{ (\tilde{C} t^2)^{i/2}}{\lfloor i / 2 \rfloor !} \leq 1 + (\sqrt{\tilde{C}} + 1)\sum_{j=1}^{\infty} \frac{(\tilde{C} t^2)^{j}}{ j !} \leq \sum_{j=0}^{\infty} \frac{(\bar{C} t^2)^{j}}{ j !} = \exp(\bar{C} t^2). \nonumber
 \end{eqnarray}
The left-hand side of (\ref{eq_926}) follows by Jensen's inequality.
We may clearly assume that $n \geq C^{\prime}$ when proving the
right-hand side inequality of (\ref{eq_926}) (otherwise, $1 - C
n^{-1/5}$ can be made negative, for an appropriate choice of
a universal constant $C$). We use (\ref{eq_933}) for the value
$$ t = \frac{L}{\sqrt{n}} = C n^{\frac{2 \alpha + 2}{5 \alpha + 20} - \frac{1}{2}} \leq C n^{-1/10} < 1, $$
to conclude that
$$ \frac{\tilde{M}(|x_0|)}{\exp(M(|x_0|))}  = \frac{\int_{SO(n)}  \exp \left( M_{f_{X+Y}, E_0, x_0}(U) \right) d\mu_n(U)}{\exp (M(|x_0|))}  $$
$$ = \int_{SO(n)}  \exp \left( M_{f_{X+Y}, E_0, x_0}(U) - M(|x_0|)
\right) d\mu_n(U) \leq \exp(\hat{C} n^{-1/5}). $$ Taking logarithms
of both sides completes the proof. \hfill $\square$

\medskip Let $X, Y, \alpha, \lambda, \ell$ be as in Lemma \ref{concent}.
We  choose a slightly different normalization. Define
\begin{equation}
Z = \frac{X+Y}{\sqrt{1 + n^{-\lambda \alpha}}},
\end{equation}
and denote by $f_Z$ the corresponding density. Clearly $f_Z$ is
isotropic and log-concave. Next we define, for $x \in E_0$,
\begin{equation} \label{eq_942_}
\tilde{M}_1(|x|) := \int_{SO(n)} \pi_{E_0} (f_{Z} \circ U)(x) d \mu_n(U).
\end{equation}
Our goal is to show that  the
following estimate holds:
\begin{equation} \label{medeq}
\left | \frac{\tilde{M}_1(|x|)}{\gamma_\ell[1](x)} - 1 \right | < C_1 n^{-c_1}
\end{equation}
for all $x \in \R^\ell$ with $|x| < c_2 n^{c_2}$ for some universal
constants $C_1,c_1,c_2>0$.

\medskip We write $S^{n-1} = \left \{ x \in
\RR^n ; |x| = 1 \right \}$, the unit sphere in $\RR^n$. Define:
\begin{equation}
\tilde{f}_Z(x) = \int_{S^{n-1}} f_Z(|x| \theta) d \sigma_n(\theta)  = \int_{SO(n)} f_Z(U x) d\mu_n(U),
\ \ \ \ \ \ \ (x \in  \RR^n)
\end{equation}
where $\sigma_n$ is the unique rotationally-invariant probability measure on $S^{n-1}$.
Since $\tilde{f}_Z$ is spherically symmetric, we shall also use the
notation $\tilde{f}_Z (|x|) = \tilde{f}_Z(x)$.
Clearly, for any $x \in E_0$,
\begin{equation} \label{sphsymm}
\tilde{M}_1(|x|) = \int_{SO(n)} \pi_{E_0} (f_{Z} \circ U)(x) d \mu_n(U)  =
\int_{SO(n)} \pi_{E_0} (\tilde{f_{Z}} \circ U)(x) d \mu_n(U) = \pi_{E_0}(\tilde{f_Z})(x).
\end{equation}
We will use the following thin-shell estimate, proved in \cite[Theorem 1.3]{K2}:
\begin{prop} \label{prop_712}
Let $n \geq 1$ be an integer and let $X$ be an isotropic random vector in $\R^n$ with
a log-concave density. Then,
\begin{equation}
\P \left \{ \left | \frac {|X|}{\sqrt n} - 1 \right | \geq \frac{1}{n^{1 / 15}} \right \} < C \exp \left( -c n^{1/15} \right)
\end{equation}
where $C,c > 0$ are universal constants.
\end{prop}
Applying the above for $f_Z$, denoting $\eps = n^{-1/15}$, and defining
$$A = \{x \in \R^n ; ~ \sqrt n(1 - \eps) \leq |x| \leq \sqrt n(1 + \eps) \}, $$ we get,
\begin{equation} \label{thinshf}
\int_A f_Z(x) dx > 1 - C e^{-c n^{1/15}}.
\end{equation}
From the definition of $\tilde{f}_Z$, it is clear that the above
inequality also holds when we replace $f_Z$ with $\tilde{f}_Z$. In
other words, if we define
\begin{equation}\label{defg}
g(t) = t^{n-1} \omega_n \tilde{f}_Z(t) \ \ \ \ \ \ \ \ \ \ \ \ \ (t \geq 0)
\end{equation}
where $\omega_n$ is the surface area of the unit sphere $S^{n-1}$ in
$\mathbb{R}^n$,
and use integration in polar coordinates, we get
\begin{equation} \label{concg}
1 \geq \int_{\sqrt n(1 - \eps)}^{\sqrt n(1 + \eps)} g(t) dt > 1 - C e^{- c n^{1/15}}.
\end{equation}
\medskip
Our next step is to apply the methods in Sodin's paper \cite{Sod} in order to prove
 a generalization of \cite[Theorem 2]{Sod},
for a multi-dimensional marginal rather then a one-dimensional marginal. Our estimate will be rather crude, but suitable for our needs. \\
Denote by $\sigma_{n, r}$ the unique rotationally-invariant probability measure on the
Euclidean sphere of radius $r$ around the origin in $\R^n$.
A standard calculation shows that the density of an
 $\ell$-dimensional marginal of $\sigma_{n,r}$ is given by the following formula:
\begin{equation}
\psi_{n,\ell,r}(x) = \psi_{n,\ell,r}(|x|) := \Gamma_{n,\ell} \frac{1}{r^\ell} \left ( 1 - \frac{|x|^2}{r^2} \right)^{\frac{n-\ell-2}{2}} 1_{[-r, r]}(|x|)
\end{equation}
where
\begin{equation}
\Gamma_{n,\ell} = \left(\frac{1}{\sqrt{\pi}} \right)^\ell \frac{\Gamma(\frac n 2)}{\Gamma(\frac{n-\ell}{2})}
\end{equation}
and where $1_{[-r,r]}$ is the characteristic function of the interval $[-r,r]$.
(see for example \cite[remark 2.10]{DF}). When $\ell << \sqrt{n}$ we have $\Gamma_{n,\ell} \left( \frac{2 \pi}{n} \right)^{\ell/2} \approx 1$. By the definition (\ref{defg}) of $g$, and since $\tilde{f}_Z$ is spherically symmetric,
we may write
\begin{equation} \label{eq_703}
\pi_{E_0}(\tilde{f}_Z)(x) = \int_0^\infty \psi_{n,\ell,r}(|x|) g(r) dr \ \  \ \ \ \ \ \ \ \ (x \in E_0).
\end{equation}
Indeed, the measure whose density is $\tilde{f}_Z$ equals $\int_0^{\infty}  g(r) \sigma_{n,r} dr$,
hence its marginal onto $E_0$ has density $x \mapsto  \int_0^\infty \psi_{n,\ell,r}(x) g(r) dr$.
We will show that the above density is approximately gaussian for $x \in E_0$ when
 $|x|$ is not too large.
But first we need the following technical lemma:
\begin{lem} \label{restrictsect}
Let $g$ be the density defined in (\ref{defg}), and suppose that $n \geq C'$
and $\ell \leq n^{1/20}$.
For $\eps = n^{-1/15}$ denote $U = \{ t > 0 ; t < (1 - \eps) \sqrt{n} \ \text{or} \ t > (1 + \eps) \sqrt{n} \}$. Then,
\begin{equation} \label{thinsh}
\int_{U} t^{-\ell} g(t) dt < C' \exp \left( - c' n^{1/15} \right).
\end{equation}
Here, $c',C'>0$ are universal constants.
\end{lem}

\textbf{Proof}: Define for convenience,
\begin{equation} \label{defh}
h(t) = t^{-\ell} g(t).
\end{equation}
Denote $$ A = \left[0, \frac 1 {n^2} \right], \ \ \ \ \ B = \left[\frac 1 {n^2}, \sqrt n(1 - \eps) \right] \cup \left[\sqrt n(1 + \eps), \infty \right), $$
and write
\begin{equation} \label{twoterms}
\int_{U} h(t) dt = \int_{A} h(t) dt + \int_{B} h(t) dt.
\end{equation}
We estimate the two terms separately. For $t>\frac 1 {n^2}$ we have
\begin{equation}
h(t) / g(t) < (n^{2 \ell}) = e^{ 2 \ell \log n}.
\end{equation}
Thus we can estimate the second term as follows:
\begin{equation} \label{lemterm2}
\int_{B} h(t) dt < e^{2 \ell \log
n} \int_{B} g(t) dt
<  e^{2 \ell \log n} C e^{-c n^{1/15}} < C e^{-\frac{1}{2} c n^{1/15}},
\end{equation}
where for the second inequality we apply the reformulation (\ref{concg}) of Proposition \ref{prop_712}
(recall that  $\eps = n^{-1/15}$ and that $\ell < n^{1/20}$).

\medskip
To estimate the first term in the right-hand side of (\ref{twoterms}), we use the fact that
$f_Z$ is isotropic and log concave, so we can use a crude bound for the
isotropic constant (see e.g. \cite[Corollary 4.3]{K3} or
\cite[Theorem 5.14(e)]{LV}) which gives
$\sup_{\R^n} f_Z < e^{n \log n}$, thus, also $\sup_{\R^n}
\tilde{f}_Z < e^{n \log n}$. Hence we can estimate
\begin{equation} \label{lemterm1}
\int_{A} h(t) dt = \int_{0}^{\frac 1 {n^2}}
t^{-\ell} g(t) dt
=  \int_{0}^{\frac 1 {n^2}} t^{n - \ell-1} \omega_n \tilde{f}_Z(t) dt
\end{equation}
\begin{displaymath}
< n^{-2(n - \ell )} \omega_n \sup \tilde{
f_Z} < e^{- 1.5 n \log n  + n \log n} < e^{-n},
\end{displaymath}
as $\omega_n < C$.
The combination of (\ref{lemterm2}) and
(\ref{lemterm1}) completes the proof. $\hfill \square$

\medskip We are now ready to show that the marginals of $\tilde{f}_Z$ are approximately gaussian.
Our desired bound (\ref{medeq}) is contained in the following lemma.

\begin{lem} \label{lem_707}
Let $1 \leq \ell \leq n$ be integers, with $n \geq C$ and $\ell \leq n^{1/20}$.
Let $g:\R^+ \to \R^+$ be a function that satisfies (\ref{concg}) and (\ref{thinsh}). Then we have,
\begin{equation}
\left| \frac{\tilde{M}_1(|x|)}{\gamma_\ell[1](x)} - 1 \right| = \left |
\frac{\int_0^\infty \psi_{n,\ell,r} (|x|) g(r) dr}{\gamma_\ell[1](x)} - 1
\right | < C n^{-1/60} \label{eq_702}
\end{equation}
for all $x \in \R^\ell$ with $|x| < 2 n^{\frac{1}{40}}$ where $C>0$ is a
universal constant.
\end{lem}
\textbf{Proof:} The left-hand side equality in (\ref{eq_702}) follows at once from (\ref{sphsymm})
and (\ref{eq_703}). We move to the proof of the right-hand side inequality.
We begin by using a well-known fact, that follows from a straightforward computation using asymptotics of
$\Gamma$-functions: for $|x| < n^{1/8}$,
\begin{equation} \label{sodlem1}
\left| \frac{\psi_{n,\ell,\sqrt n}(|x|)}{\gamma_\ell [1](x)} - 1 \right| = \left| \left( \frac{2 \pi}{n} \right)^{\ell/2} \Gamma_{n,\ell} \frac{\left(1 - \frac{|x|^2}{n} \right)^{(n - \ell-2) / 2}}{e^{-|x|^2/2}} - 1 \right|
\leq \frac{C}{\sqrt n}
\end{equation}
(We omit the details of the simple computation. An almost identical computation is done, for example, in \cite[Lemma 1]{Sod}.
Note that in addition to the computation there, we have to use, e.g.,
 Stirling's formula to estimate the constants $\eps_n$).
Using the above fact (\ref{sodlem1}), we see that it suffices to prove the following
inequality:
\begin{equation} \label{reducttosod}
\left | \frac{\int_0^\infty \psi_{n,\ell,r} (|x|) g(r)
dr}{\psi_{n,\ell,\sqrt n}(|x|) } - 1 \right | < C n^{-\frac{1}{60}}
\end{equation}
for all $x \in \R^\ell$ with $|x| < 2 n^{\frac{1}{40}}$.
To that end, fix $x_0 \in \R^\ell$ with $|x_0| < 2 n^{\frac{1}{40}}$,
define $$ A = [\sqrt n (1 - n^{-\frac{1}{15}}), \sqrt n (1 +
n^{-\frac{1}{15}})], \ \ \ \ \ \ B=[0, \infty) \setminus A, $$ and
write
\begin{equation} \label{tterms}
\int_0^\infty \psi_{n,\ell,r} (|x_0|) g(r) dr  = \int_A \psi_{n,\ell,r} (|x_0|) g(r) dr + \int_{B} \psi_{n,\ell,r} (|x_0|) g(r) dr.
\end{equation}
We estimate the two terms separately. For the second term, we have,
\begin{eqnarray}
\nonumber
\lefteqn{
\int_{B} \psi_{n,\ell,r} (|x_0|) g(r) dr = \Gamma_{n,\ell} \int_{B} \frac{1}{r^\ell} \left ( 1 - \frac{|x_0|^2}{r^2} \right)^{\frac{n-\ell-2}{2}} 1_{[-r, r]}(|x_0|) g(r) dr} \\ &
< & \Gamma_{n,\ell} \int_{B} \frac{1}{r^\ell} g(r) dr < \Gamma_{n,\ell} C e^{-c n^{1/15}},
\phantom{aaaaaaaaaaaaaaaaaaaaaaaaaaaaa}
\end{eqnarray}
where the last inequality follows from (\ref{thinsh}).
Therefore,
\begin{eqnarray}\label{inta1}
\lefteqn{ \frac{\int_{B} \psi_{n,\ell,r} (|x_0|) g(r) dr}{\psi_{n,\ell,\sqrt
n}(|x_0|)} <  \frac{C e^{-c n^{1/15}}}{ (\frac{1}{\sqrt
n})^\ell \left ( 1 - \frac{|x_0|^2}{n} \right)^{\frac{n-l-2}{2}} } } \\ & < &
C e^{-c n^{1/15} + |x_0|^2 + \frac{1}{2} \ell \log n} < C e^{-n^{1/20}}. \nonumber
\end{eqnarray}
To estimate the first term on the right-hand side of (\ref{tterms}), we will show that the following inequality holds:
\begin{equation}
\left |\frac{\int_A \psi_{n,\ell,r} (|x_0|) g(r) dr}{\psi_{n,\ell, \sqrt
n}(|x_0|)} - 1 \right | < C n^{-1/60}
\end{equation}
for some constant $C > 0$. For $r>0$ such that $\frac{|x_0|^2}{r^2} < \frac{1}{2}$, we have,
\begin{equation}
\left |\frac{d}{dr} \log \psi_{n,\ell,r} (|x_0|) \right | = \left |- \frac{\ell}{r} + (n-\ell-2) \frac{|x_0|^2}{r^3} \frac{1}{\left ( 1 - \frac{|x_0|^2}{r^2} \right)} \right | < \frac{\ell}{r} + 2 n \frac{|x_0|^2}{r^3}.
\end{equation}
Recalling that $|x_0| < 2 n^{\frac{1}{40}}$ and $\ell \leq n^{1/20}$, the above estimate gives
that for all $r \in [\frac{1}{2} \sqrt n, \frac{3}{2} \sqrt n]$,
\begin{equation}
\left | \frac{d}{dr} \log \psi_{n,\ell,r} (|x_0|) \right | <
2n^{\frac{1}{20} - \frac{1}{2}} +  16 n^{1 + \frac{1}{20} - \frac{3}{2}}
< C n^{-\frac{9}{20}}
\end{equation}
which gives, for $r \in [\frac{1}{2} \sqrt n, \frac{3}{2} \sqrt n]$,
\begin{equation}
\left | \frac{\psi_{n,\ell,r}(|x_0|)}{\psi_{n,\ell,\sqrt n}(|x_0|)} - 1 \right | < C n^{-\frac{9}{20}}|r - \sqrt n| .
\end{equation}
Recall that for $r \in A$ we have $|r - \sqrt n| \leq n^{\frac{13}{30}}$. Hence the last estimate yields,
\begin{equation}
\left |\frac{\int_A \psi_{n,\ell,r} (|x_0|) g(r) dr}{\psi_{n,\ell, \sqrt
n}(|x_0|) \int_A g(r) dr} -  1 \right | < C n^{-\frac{9}{20}}
n^{\frac{13}{30}} = C n^{-\frac{1}{60}}.
\end{equation}
Combining the last inequality with (\ref{concg}), we get
\begin{equation} \label{inta2}
\left |\frac{\int_A \psi_{n,\ell,r} (|x_0|) g(r) dr}{\psi_{n,\ell, \sqrt n}(|x_0|)} -  1 \right | <
\tilde{C} e^{-c n^{\frac{1}{15}}} + C n^{-\frac{1}{60}} < C'
n^{-\frac{1}{60}}.
\end{equation}
From (\ref{inta1}) and (\ref{inta2}) we deduce (\ref{reducttosod}), and the lemma is proved. $\hfill \square$

\medskip
Recall the definitions (\ref{eq_942}) and (\ref{eq_942_}) of $\tilde{M}(|x|)$ and $\tilde{M}_1(|x|)$; the only difference is the normalization of $X+Y$.
By an easy scaling argument, we deduce from Lemma \ref{lem_707} that when $n \geq C$,
\begin{equation} \label{finalmt}
\left | \frac{\tilde{M}(|x|)}{\gamma_\ell[1 + n^{-\lambda \alpha}](x)} - 1
\right | < C_1 n^{-\frac{1}{60}}
\end{equation}
for all $x \in \R^\ell$ with $|x| < n^{\frac{1}{40}}$, for $C_1 > 0$
a universal constant.  By plugging (\ref{eq_926}) and (\ref{finalmt})
into Lemma \ref{concent},
we conclude the following:

\begin{prop} \label{convgauss}
Let $1 \leq \ell \leq n$ be integers. Let $0 < \alpha < 10^5$ and denote $\lambda = \frac{1}{5 \alpha + 20}$.
Assume that $\ell \leq n^\lambda$. Suppose that
$f:\mathbb{R}^n \to [0,\infty)$ is a log-concave function that is the density of an
isotropic random vector.  Define $g = f
* \gamma_n[n^{-\lambda \alpha}]$, the convolution of $f$ and $\gamma_n[n^{-\lambda \alpha}]$.
Let $E \in G_{n,\ell}$ be a random
subspace. Then, with probability greater than $1-C
e^{-c n^{1/10}}$ of selecting $E$, we have
\begin{equation}
\left | \frac{\pi_E(g)(x)}{\gamma_\ell [1 + n^{-\lambda \alpha}](x)} - 1
\right | \leq C n^{-\lambda}
\end{equation}
for all $x \in E$ with $|x| < n^{\lambda / 2}$, where $C > 0$ is a universal constant.
\end{prop}

We did not have to explicitly assume that $n \geq C$ in
Proposition \ref{convgauss}, since otherwise the proposition is vacuously true.
In the next section we will show that the above estimate still holds
without taking the convolution, perhaps with slightly worse
constants.

\section{Deconvolving the Gaussian}
\label{sec3}
Our goal in this section is to establish the following principle:
Suppose that $X$ is a random vector with a log-concave density,
and that $Y$ is an independent, gaussian random vector whose covariance
matrix is small enough with respect to that of $X$.
 Then, in the case where $X + Y$ is approximately gaussian,
the density of $X$ is also approximately gaussian,
in a rather large domain. We begin with a lower bound for the density of $X$. \\
(Note that the notation $n$ in this section
corresponds to the dimension of the subspace, that was denoted by $\ell$ in the previous  section.)
\begin{lem} \label{lower}
Let $n \geq 1$ be a dimension, and let $\alpha,\beta, \eps, R
> 0$. Suppose that $X$ is an isotropic random vector in $\RR^n$ with
a log-concave density, and that $Y$ is an independent gaussian
random vector in $\RR^n$ with mean zero and covariance matrix $\alpha Id$.
Denote by $f_X$ and $f_{X+Y}$ the respective densities. Suppose
that,
\begin{equation} \label{convbound}
f_{X+Y}(x) \geq (1 - \eps) \gamma_n[1 + \alpha](x)
\end{equation}
for all $|x| \leq R$. Assume that $\alpha \leq c_0 n^{-8}$ and that
\begin{equation} \label{xicond}
100 (2n)^{\max \{ 3 \beta, 3 / 2 \}} \alpha^{1/4} < \eps < \frac{1}{100}.
\end{equation}
Then,
\begin{equation}
f_{X}(x) \geq (1 - 6 \eps) \gamma_n[1](x)
\end{equation}
for all $x \in \RR^n$ with $|x| \leq \min \left \{R - 1, (2n)^\beta \right \}$. Here, $0 < c_0 < 1$
is a universal constant.
\end{lem}

\textbf{Proof:} Suppose first that $f_X$ is positive everywhere in $\RR^n$.
Fix $x_0 \in \RR^n$ with $|x_0| \leq \min \{ R - 1, (2n)^{\beta} \}$.
Assume that $\eps_0 > 0$
is such that
\begin{equation} \label{defep0}
 f_X (x_0) < (1 - \eps_0) \gamma_n[1](x_0).
\end{equation}
To prove the lemma (for the case where $f_X$ is positive everywhere) it suffices to show that
\begin{equation}
\eps_0 \leq 6 \eps.
\label{goal}
\end{equation}
Consider the level set $L = \{x \in \RR^n;  f_X(x) \geq f_X(x_0)\}$.
Then $L$ is convex and bounded, as $f_X$ is log-concave and integrable (here we used the
fact that $f_X(x_0) > 0$).
Let $H$ be an affine hyperplane that supports $L$ at its boundary point $x_0$,
and denote by $D$ the open ball of radius $\alpha^{1 / 4}$ tangent to $H$ at $x_0$,
that is disjoint from the level set $L$.
By definition, $f_X(x) < f_X(x_0)$ for $x \in D$.
Denote the center of $D$ by $x_1$. Then, $\left | x_1 - x_0 \right | \leq \alpha^{1 / 4}$
with $|x_0| \leq (2n)^\beta$,
and a straightforward computation yields
\begin{equation} \label{norm1}
\left||x_1|^2 - |x_0|^2 \right | \leq  \left( 2 (2n)^{\beta} + \alpha^{1 / 4} \right) \alpha^{1/4} \leq  \frac{\eps}{2},
\end{equation}
where we used (\ref{xicond}).
Note that $|x_1| \leq |x_0| + \alpha^{1/4} \leq R$.
Apply the last inequality and (\ref{convbound}) to obtain,
\begin{equation} \label{main1}
f_{X+Y}(x_1) \geq (1 - \eps) \gamma_n[1+\alpha](x_0) e^{\frac{|x_0|^2 - |x_1|^2}{2(1+\alpha)}} > (1 - 2 \eps) \gamma_n[1+\alpha](x_0).
\end{equation}
By definition,
\begin{equation} \label{eq_958}
f_{X+Y}(x_1) = \int_{\R^n} f_X(x) \gamma_n [\alpha] (x_1-x) dx =
\end{equation}
\begin{displaymath}
\int_{x \in D} f_X(x) \gamma_n [\alpha] (x_1-x) dx + \int_{x \notin D} f_X(x) \gamma_n [\alpha]
(x_1-x) dx.
\end{displaymath}
We will estimate both integrals. First, recall that $f_X(x) < f_X(x_0)$ for $x \in D$ and use (\ref{defep0})
to deduce
\begin{equation} \label{int1}
\int_{x \in D} f_X(x) \gamma_n [\alpha] (x_1-x) dx < f_X(x_0) < (1 - \eps_0) \gamma_n[1](x_0).
\end{equation}
For the integral outside $D$ a rather rough estimate would suffice. We may write,
\begin{equation}
\int_{x \notin D} f_X(x) \gamma_n [\alpha] (x_1 -x) dx < \mathbb{P} \left(|G_n| \geq \frac{1}{\alpha^{1/4}} \right) \sup_{\R^n} f_X
\end{equation}
where $G_n \sim \gamma_n[1]$ is a standard gaussian random vector.
To bound the right-hand side term, we shall use a standard
tail bound for the norm of a gaussian random vector,
\begin{equation} \label{gausstail}
\P(|G_n| > t \sqrt n) < C e^{-ct^2},
\end{equation}
and the following crude bound for the isotropic constant of $f_X$ (see e.g \cite[Theorem 5.14(e)]{LV}),
\begin{equation} \label{lfbound}
\sup_{\R^n} f_X < e^{ \frac{1}{2} n \log n + 6 n} < e^{C n \log n}.
\end{equation}
Consequently,
\begin{equation}\label{int2}
\int_{x \notin D} f_X(x) \gamma_n [\alpha] (x_1-x) dx < Ce^{-c n^{-1} \alpha^{-1/2} } e^{C n \log n} < e^{- \alpha^{- 1/ 3}},
\end{equation}
for an appropriate choice of a sufficiently small universal
constant $c_0 > 0$ (so that all other constants are absorbed).
Combining (\ref{eq_958}), (\ref{int1}) and (\ref{int2}) gives
\begin{equation} \label{eq_145}
f_{X+Y}(x_1) < \left(1 - \eps_0 + \frac{e^{-\alpha^{-1/3}}}{\gamma_n[1](x_0)} \right) \gamma_n[1](x_0).
\end{equation}
Using the fact that $n + (2n)^{2 \beta} < \frac{\alpha^{-1/3}}{2}$, which follows
easily from our assumptions, we have
\begin{equation}
\frac{e^{-\alpha^{-1/3}}}{\gamma_n[1](x_0)} =
e^{\frac{|x_0|^2}{2} + \frac{n}{2} \log (2 \pi) - \alpha^{-1/3}} < e^{-\frac{1}{2} \alpha^{-1/3}} \leq 2 \alpha^{1 / 3} < \frac{\eps}{2} < \frac{\eps_0}{2} \label{eq_146}
\end{equation}
(for the last inequality, note that if $\eps_0 < 6 \eps$ then we
have nothing to prove. So we
can assume that $\eps_0 > \eps$). From (\ref{eq_145}) and (\ref{eq_146}) we obtain
the bound
\begin{equation} \label{main2}
f_{X+Y}(x_1) <  \left(1 - \frac{\eps_0}{2} \right) \gamma_n[1](x_0).
\end{equation}
Combining (\ref{main1}) and (\ref{main2}) we get,
\begin{equation}
(1 - 2 \eps) \gamma_n[1+\alpha](x_0) < \left(1 - \frac{\eps_0}{2} \right) \gamma_n[1](x_0).
\end{equation}
A calculation yields,
\begin{equation} \label{difvar}
\frac{\gamma_n[1](x_0)}{\gamma_n[1+\alpha](x_0)} \leq \frac{\gamma_n[1](0)}{\gamma_n[1+\alpha](0)} = (1 + \alpha)^{\frac{n}{2}} < 1 + \eps.
\end{equation}
From the above two inequalities, we finally deduce,
\begin{equation}
\frac{1 - \eps_0 / 2}{1 - 2 \eps} > \frac{1}{1 + \eps} > 1 - \eps \ \ \ \ \ \ \Rightarrow \ \ \ \ \ \ \eps_0 < 6 \eps,
\end{equation}
which proves (\ref{goal}). The lemma is proved, under the additional assumption that $f_X$
never vanishes. The general case follows by a standard approximation argument.
 $\hfill \square$

\medskip
After proving a lower bound, we move to the upper bound. We will show that if we add to the
requirements of the previous lemma an assumption that the density of $f_{X+Y}$
is bounded from above, then we can provide an upper bound for $f_X$.
\begin{lem} \label{upper}
Let $n,X,Y, \alpha, \beta, \eps, R, c_0$ be defined as in Lemma \ref{lower}, and suppose that all
the conditions of Lemma \ref{lower} are satisfied. Suppose that in addition, we have the following upper bound for $f_{X+Y}$:
\begin{equation}\label{co1}
f_{X+Y}(x) < (1 + \eps) \gamma_n[1 + \alpha](x)
\end{equation}
for all $|x| < R$. Then we have:
\begin{equation}
f_{X}(x) < (1 + 8 \eps) \gamma_n[1](x)
\end{equation}
for all $x$ with $|x| < \min \{ (2n)^{\beta}, R \} - 3$.
\end{lem}
\textbf{Proof:}
Denote $F(x) = -\log f_X(x)$. Again we use the upper bound for the supremum of the density (\ref{lfbound}),
\begin{equation} \label{upperf}
F(x) > 6 n - \frac{1}{2} n \log n  > -  n \log n, ~~~\forall x \in \R^n.
\end{equation}
Use the conclusion of Lemma \ref{lower} to deduce that for $|x| < \min \{(2n)^{\beta}  , R\} -
1$ the following holds:
\begin{equation} \label{lowerf}
F(x) < - \log \left( \frac 1 2 \gamma_n[1](x) \right) < \log 2 + \frac n 2 \log (2\pi) + (2n)^{2 \beta} < 3 (2n)^{\max\{2 \beta, \frac{3}{2}\}}.
\end{equation}
 Next we will show that for $x, y
\in A = \left \{x  \in \RR^n ; ~ |x| < \min \{(2n)^{\beta}, R\} -
2~ \right \}$, the following Lipschitz condition holds:
\begin{equation} \label{lipf}
|F(x) - F(y)| \leq 5 (2n)^{\max\{2 \beta, \frac{3}{2}\}} |x - y|.
\end{equation}
To that end, denote $a = 5 (2n)^{\max\{2 \beta, \frac{3}{2}\}}$ and suppose by
contradiction that $x,y \in A$ are such that
\begin{equation}
F(y) - F(x) > a |y - x|.
\end{equation}
Since $F(y) - F(x) < a$ (as implied by (\ref{upperf}) and
(\ref{lowerf})), we have $|y - x| < 1$ and for the point
$$ y_1 := x + \frac{y - x}{|y-x|}, $$
we have, using the convexity of $F$,
$$F(y_1) - F(x) \geq \frac{F(y) - F(x)}{|y-x|} >  a.$$
Note that $|y_1| \leq |x| + 1 < \min \{(2n)^{\beta}, R\} -
1$, so we
get a contradiction to (\ref{upperf}) and (\ref{lowerf}).
This proves (\ref{lipf}). \\
Therefore, given two points $x, x_0 \in A$ such that $|x_0 - x| <
\alpha^{1/ 4}$, (\ref{lipf}) implies,
\begin{equation}
\left| F(x_0) - F(x) \right| < 5 \alpha^{1 / 4} (2n)^{\max\{2 \beta,
3 /2 \}} < \eps / 20.
\end{equation}
Recall that $F = -\log f_X$, hence the above translates to
\begin{equation} \label{lipconc}
|f_X(x_0) - f_X(x)| < 2 (e^{\eps/20} - 1) f_X(x_0) < \frac{\eps}{4}
f_X(x_0).
\end{equation}
Now, suppose $x_0 \in \R^n$ and $0 < \eps_0 < 1$ are such that
\begin{equation} \label{asumep0}
f_X(x_0) > (1 + \eps_0) \gamma_n [1](x_0),
\end{equation}
with $|x_0| < \min \{R, (2n)^\beta\} - 3$. Again, to prove the Lemma it suffices
to show that in fact $\eps_0 < 8 \eps$.
Let $D$ be a ball of radius $\alpha^{1/4}$ around $x_0$. \\
Since we can assume that $\eps_0 > \eps$ (otherwise, there is nothing
to prove), we deduce from
(\ref{lipconc}) and (\ref{asumep0}) that for all $x \in D$,
\begin{equation}
f_X(x) > \left( 1 - \frac{\eps_0}{4} \right) \left(1 + \eps_0
\right) \gamma_n [1](x_0) > \left(1 + \frac{\eps_0}{2} \right)
\gamma_n [1](x_0).
\end{equation}
Thus,
\begin{eqnarray} \label{co2}
 \lefteqn{f_{X+Y}(x_0) = \int_{\RR^n} f_X(x) \gamma_n [\alpha] (x_0-x) dx } \\
\nonumber & > & \int_{x \in D} f_X(x) \gamma_n [\alpha] (x_0-x) dx \\
& > & \left(1 + \frac {\eps_0}{2} \right) \gamma_n [1](x_0) \cdot \left(1 - \mathbb{P} \left(|G_n| > \frac{1}{\alpha^{1/4}} \right) \right) > \left(1 + \frac {\eps_0}{3} \right)\gamma_n [1](x_0), \nonumber
\end{eqnarray}
where in the last inequality we used the estimate (\ref{gausstail}) and
the assumption $\eps_0 > \eps$. Now, a computation yields,
\begin{equation} \label{difvar2}
\frac{\gamma_n[1+\alpha](x_0)}{\gamma_n[1](x_0)} < e^{\frac{1}{2}( |x_0|^2 - \frac{|x_0|^2}{1 + \alpha})}  = e^{\frac{1}{2} |x_0|^2 \frac{\alpha}{1 + \alpha}} < e^{(2n)^{2 \beta} \alpha} < 1 + \eps.
\end{equation}
We thus obtain, combining (\ref{co1}) and (\ref{co2}) and using (\ref{difvar2}), that
\begin{displaymath}
\frac{1 + \eps_0 / 3}{1 + \eps} < \frac{\gamma_n [1 + \alpha](x_0)}{\gamma_n [1](x_0)} < 1 + \eps,
\end{displaymath}
so $\eps_0 < 8 \eps$, and the proof of the lemma is complete. $\hfill \square$

\medskip
The combination of the two above lemmas gives us the desired estimate for the density of $X$,
as advertised in the beginning of this section.

\section{Proof of main theorem}

\textbf{Proof of Theorem \ref{thm_main}:}
We may clearly assume that $n$ exceeds some positive universal constant (otherwise,
take $\cE = \emptyset$).
Let $1 \leq \ell \leq n^{1/100}$ be an integer, and let $\delta \geq 0$ be
such that $\ell=n^{\delta}$.
Set $\alpha = 10$ and $\lambda = \frac{1}{5 \alpha + 20} = \frac{1}{70}$.
Let $Y$ to be a gaussian random vector in $\RR^n$ with mean zero and  covariance
matrix $n^{-\alpha \lambda} Id$,  independent of $X$.
We first apply Proposition \ref{convgauss}
for the random vector $X+Y$ with parameters $\ell$ and $\alpha$
(noting that $\ell \leq n^{1/100} \leq n^{\lambda}$).
According to the conclusion of that proposition,
if $E$ is a random subspace of dimension $\ell$, then
\begin{equation}
\left | \frac{\pi_E (f_{X+Y})(x)}{\gamma_n [1 + n^{- \alpha \lambda}] (x)} - 1
\right | \leq C n^{-1/100}, \label{eq_1016}
\end{equation}
for all $x \in E$ with $|x|< n^{\frac 1 {200}}$, with probability
greater than $1 - Ce^{-c n^{1/10}}$ of choosing $E$.

\medskip Next, we apply Lemma \ref{lower} and Lemma \ref{upper} in the $\ell$-dimensional subspace $E$,
with the parameters $\alpha = n^{-10\lambda} \leq n^{-1/20} \ell^{-8}$,
$\beta=\frac{1}{600 (\delta + 1/\log_2 n)}$, $R = n^{1/200}$, $\eps =  C n^{-1/100}$
where $C$ is the constant from (\ref{eq_1016}).
It is straightforward to verify that the requirements of these two lemmas hold,
since $n$ may be assumed to exceed a given universal constant.
According to the conclusions of Lemma \ref{lower} and Lemma \ref{upper},
for any $x \in E$ with $|x|< n^{\frac 1 {700}}$,
$$
\left | \frac{\pi_E (f_{X})(x)}{\gamma_n [1] (x)} - 1
\right | \leq C' n^{-1/100}. $$
This completes the proof. \hfill $\square$

\medskip {\bf Remark.} The numerical values of the exponents
$c_1,c_2,c_3,c_4$ provided by our proof
of Theorem \ref{thm_main} are far from optimal. The theorem is tight only in
the sense that the power-law dependencies on $n$ cannot be improved to, say, exponential dependence.
The only constant among $c_1,c_2,c_3,c_4$ for which the best value is essentially known to us is $c_2$.
It is clear from the proof that $c_2$ can be made arbitrarily close to $1$ at the expense
of decreasing the other constants. Note also that necessarily  $c_4 \leq 1/4$,
as is shown by the example where $X$ is distributed uniformly in a Euclidean ball
(see \cite[Section 4.1]{Sod}).

\bigskip {\noindent School of Mathematical Sciences, Tel-Aviv University, Tel-Aviv
69978, Israel \\  {\it e-mail address:}
\verb"roneneldan@gmail.com" }

\bigskip
{\noindent Department of Mathematics, Princeton
university, Princeton, NJ 08544, USA \\ {\it e-mail address:}
\verb"bklartag@princeton.edu"

\end{document}